\documentclass[10pt]{amsart}

\usepackage{amssymb,amsmath,latexsym}
\usepackage{charter,eucal}
\usepackage{amssymb}
\usepackage{amscd}
\usepackage{tikz}
\usetikzlibrary{positioning}
\usetikzlibrary{arrows}
\usepackage{hyperref,mathrsfs}

\newtheorem{theorem}{Theorem }[section]

\newtheorem{lemma}[theorem]{Lemma}

\newtheorem{remark}[theorem]{Remark}
\newtheorem{corollary}[theorem]{Corollary}
\newtheorem{proposition}[theorem]{Proposition}
\newtheorem{observation}[theorem]{Observation}
\newtheorem{conjecture}[theorem]{Conjecture}

\newcommand{\bt}{\begin{theorem}}
\newcommand{\et}{\end{theorem}}
\newcommand{\bc}{\begin{corollary}}
\newcommand{\bl}{\begin{lemma}}
\newcommand{\ec}{\end{corollary}}
\newcommand{\el}{\end{lemma}}
\newcommand{\bo}{\begin{observation}}
\newcommand{\eo}{\end{observation}}
\newcommand{\bp}{\begin{proposition}}
\newcommand{\ep}{\end{proposition}}
\newcommand{\br}{\begin{remark}}
\newcommand{\er}{\end{remark}}
\newcommand{\brt}{\begin{result}}
\newcommand{\ert}{\end{result}}
\newcommand{\eop}{\hspace*{\fill}{\footnotesize $\blacksquare$}}

\def\I{\texttt{I}}
\def\nI{\not\mathbf{I}}

\def\PG{\mathbf{PG}}

\def\Aut{\mathrm{Aut}}
\def\endo{\mathrm{end}}
\def\endot{\mathrm{end}^{\mathrm{t}}}
\def\id{\mathrm{id}}

\def\Z{\mathbb{Z}}

\newcommand{\hP}{\mathbf{P}}
\newcommand{\proj}{\mathrm{proj}}

\newcommand{\mU}{\mathcal{U}}
\newcommand{\mF}{\mathbf{F}}
\newcommand{\hF}{\mathbf{F}}

\newcommand{\he}{\mathrm{id}}
\newcommand{\mP}{\mathcal{P}}
\newcommand{\mB}{\mathcal{B}}
\newcommand{\mL}{\mathcal{L}}
\newcommand{\mQ}{\mathcal{Q}}
\newcommand{\mN}{\mathcal{N}}
\newcommand{\mO}{\mathcal{O}}
\newcommand{\mA}{\mathcal{A}}
\newcommand{\mG}{\mathcal{G}}

\newcommand{\mS}{\mathcal{S}}

\newcommand{\F}{\mathcal{F}}
\newcommand{\bF}{\mathbb{F}}

\newcommand{\mC}{\mathcal{C}}
\newcommand{\bA}{\mathbb{A}}
\newcommand{\bP}{\mathbb{P}}
\newcommand{\bT}{\mathbb{T}}
\newcommand{\N}{\mathbb{N}}

\newcommand{\vo}{\overline{0}}

\title{A translation generalized quadrangle in characteristic $\ne 0$ is linear}



\author{Koen Thas}

\thanks{}

\address{Ghent University, Department of Mathematics\\
Krijgslaan 281, S25, B-9000 Ghent, Belgium}

\email{koen.thas@gmail.com}

\date{}

\begin{document}

\maketitle

\begin{abstract}
It is a long-standing conjecture from the 1970s that every translation generalized quadrangle is {\em linear}, that is, 
has an endomorphism ring which is a division ring (or, in geometric terms, that has a {\em projective representation}). 
We show that any translation generalized quadrangle $\Gamma$ is ideally embedded in a translation quadrangle which is linear. This allows us to weakly represent any such $\Gamma$ in projective space, and moreover, to have a well-defined notion of ``characteristic'' for these objects. 
We then show that each translation quadrangle in positive characteristic indeed is linear.
\end{abstract}

\medskip
\tableofcontents

\section{Introduction}

If $\rho$ is a hyperplane of a projective space $\bP^m(\ell)$ defined over the field $\ell$, with $m \geq 2$, then the group of translations fixes $\rho$ pointwise, and acts sharply transitively on the set of points not in $\rho$. If $\Gamma$ is an axiomatic projective plane, and $\rho$ is a line for which such a translation group $A$ exists, it is necessarily unique, and abelian. A plane with such a group is by definition a {\em translation plane} \cite{Jo}.
Fix any point $u \in \Gamma \setminus \rho$, and for any line $Y$ incident with $u$, let $A_U$ be the stabilizer of $U$ in $A$. Define $\endot(A,A)$ to be the set of endomorphisms $\gamma: A \mapsto A$ such that $\gamma(A_X) \leq A_X$ for every line $X$ on $u$; then an old and basic result says that $\endot(A,A)$ is a division ring (where addition and multiplication are inherited from the endomorphism ring of $A$). This leads to the fact that one can represent $\Gamma$ in a projective space over $\endot(A,A)$, since the groups $A$, $A_X$ naturally are vector spaces over $\endot(A,A)$. This representation is the  so-called  ``projective'' or ``Bruck-Bose'' representation, and is arguably the most fundamental tool to study and construct translation planes. It was founded in Bruck and Bose their 1964 paper \cite{BrBo}.\\


Projective planes are examples of Tits buildings of rank $2$ | also called ``generalized polygons'' or ``generalized $n$-gons,'' when the gonality of the geometry is specified. Here $n \geq 3$, and the case $n = 3$ corresponds to projective planes. The notion of translation plane naturally generalizes to that of translation generalized $n$-gons | see \cite[4.9.8]{POL}, but it is easily shown that the very definition forces the gonality $n$ to be in the set $\{3,4\}$, see \cite[4.9.8]{POL}. For generalized $4$-gons, or {\em generalized quadrangles}, a very rich theory exists, similar to the case of translation planes. There also, in the same way as for translation planes, one defines an endomorphism ring $\endot(A,A) =: k$ | called ``kernel,'' and the main conjecture is whether it is a division ring so as to have projective representations at hand. (Precise definitions will be given in \S\S \ref{TGQintro}.) It is known that $k$ is a (not necessarily commutative) integral domain \cite{Jos}.

So the question is, very roughly, {\em whether all translation generalized quadrangles} (TGQs) {\em can be appropriately embedded in projective space over a division ring.} This question is the main obstacle in the current theory of  general translation quadrangles.
The origins of the problem can be traced back to the 1970s, when the conjecture was established for finite TGQs | see \cite{TGQold}, and the monographs \cite[chapter 8]{PT2} and \cite{TGQ}.

\begin{conjecture}[Linearity for TGQs]
\label{conjlinear}
Any translation generalized quadrangle has a projective representation, that is, the kernel
$k$ is a division ring.
\end{conjecture}


Throughout this paper, we will call TGQs for which the kernel is a division ring {\em linear}, as in \cite{Jos}.

\subsection{The present paper}

We say that a TGQ $(\mS,A)$ has a {\em weak projective representation} if there exists a linear TGQ $(\overline{\mS},\overline{A})$ which ideally contains 
$(\mS,A)$. In that case, since $(\overline{\mS},\overline{A})$ has a projective representation, $(\mS,A)$ can also be represented in projective space.

Assuming by way of contradiction that the kernel $k$ of $\Gamma$ is not a division ring, we will build inductively an ``ambient translation generalized quadrangle'' $\widetilde{\Gamma}$ for which the kernel $\widetilde{k}$ {\em is} a division ring, and which contains $\Gamma$ as an ideal subquadrangle. 
Our first  main result  then reads as follows.

\begin{theorem}[Weak projective representation of TGQs]
\label{conjrough}
Any translation generalized quadrangle has a weak projective representation. 
\end{theorem}

\br{\rm
Theorem \ref{conjrough} is a direct corollary of Conjecture \ref{conjlinear} if the latter were true.}
\er

The reason that $\widetilde{\Gamma}$ exists rests,
essentially, upon the following simple fact: if $G$ is a group and $\zeta: G \rightarrow G$ is an injective endomorphism, then the direct limit of the directed system
\begin{equation}
\cdots \xrightarrow{\zeta} G \xrightarrow{\zeta} G  \xrightarrow{\zeta} \cdots
\end{equation}
is again a group for which $\zeta$ is stable (that is, becomes an automorphism). [The direct limit can be naturally seen as a union of isomorphic copies
of $G$ which form a flag, and $\zeta$ acts on this limit, by the action on each of the copies.] We will introduce a similar localization process for TGQs in order to eventually produce ambient TGQs after at most a countable number of repetitions of this procedure.

For TGQs that have a weak projective representation, we can define a {\em characteristic}, namely the characteristic of the kernel of any ambient TGQ (it is independent of the choice of ambient TGQ). Our second main result then is:

\bt
A TGQ of positive characteristic is linear.
\et

We will also severely restrict the possibilities for counter examples in characteristic $0$.

\subsection{The known cases}

For some special cases Conjecture \ref{conjlinear} has been proved to be true; {\em planar} TGQs \cite{Jos2}, TGQs with a {\em strongly regular translation center} \cite{Jos2} and of course {\em finite} TGQs \cite{PT2} all satisfy the conjecture.
Also, it has been shown in \cite{Jos3} that the more restricted ``topological kernel'' of a compact connected topological TGQ is a division ring.  Due to these rather restricted partial results, it is usually \underline{assumed} that the kernel of a TGQ be a division ring.

\medskip
\br
{\rm
In \cite{Jos} it is claimed (in Corollary 3.11) that the kernel of a TGQ always is a division ring, the result which is the main purpose of the present paper. In the proof however, the author uses Proposition 3.10 which states 
that any three distinct lines on the translation point, together with any affine point, generate a ``plane-like''
subGQ. This result is not true (even not in the finite case) | in fact, only a very restricted class of TGQs
has this property. (Still, Theorem 3.11 of {\em loc. cit.} shows that if a TGQ {\em does} satisfy this property,
it indeed {\em is} linear.) The paper \cite{Jos} contains many other interesting results on infinite TGQs.\\
}
\er

\section{Some definitions}

We introduce some frequently used notions in the following paragraphs.

\subsection{Generalized quadrangles and subquadrangles}

A {\em generalized $4$-gon} or {\em generalized quadrangle} (GQ) \cite{PT2} is a point-line incidence geometry $\Gamma = (\mP,\mB,\I)$ for which  the following axioms are satisfied:

\begin{itemize}
\item[(i)]
$\Gamma$ contains no ordinary $k$-gon (as an induced subgeometry), for $2 \leq k < 4$;
\item[(ii)]
any two elements $x,y \in \mP \cup \mB$ are contained in some ordinary $4$-gon in
$\Gamma$;
\item[(iii)]
there exists an ordinary $5$-gon in $\Gamma$.\\
\end{itemize}

Note that points and lines play the same role in this definition. This is the principle of ``duality.''

Incidence geometries which satisfy (i) and (ii) but not (iii) are {\em thin} GQs.

If $x$ and $y$ are collinear points in $\Gamma$, we write $x \sim y$; in particular, $x \sim x$. 
We use the same notation for lines. 

Generalized quadrangles were introduced by Tits in his triality paper
 \cite{Ti} in order to better understand  the Chevalley groups of rank $2$.

\subsection{Order}

Each thick GQ has the property that there exist values $u > 1$ and $v > 1$ such that each line is incident with $u + 1$ points and 
each point is incident with $v + 1$ lines; see \cite[\S 1.5.3]{POL}. We say that $(u,v)$ is the {\em order} of the quadrangle, and $u, v$ are its {\em parameters}. In this paper, almost without exception the 
parameters will be infinite.

\subsection{Subquadrangles}

A {\em subquadrangle} (subGQ) of a generalized quadrangle $\Gamma = (\mP,\mB,\I)$ is a generalized quadrangle $\Gamma' = (\mP',\mB',\I')$ such that
$\mP' \subseteq \mP$, $\mB' \subseteq \mB$, and $\I' = \I \cap \Big((\mP' \times \mB') \cup (\mB '\times \mP')\Big)$. We write $\Gamma' \preceq \Gamma$ to express that $\Gamma'$ is a subGQ of $\Gamma$. If $\Gamma' \preceq \Gamma$ and if $\Gamma'$ has the property that for every point $x$ of $\Gamma'$, $U \I x$ in $\Gamma$ implies that $U$ is also a line of $\Gamma'$, then $\Gamma'$ is an {\em ideal} subquadrangle.

\subsection{Regularity}
\label{reg}



Let $x, y$ be distinct noncollinear points in a thick generalized quadrangle. Then we say that the pair $\{x,y\}$ is {\em regular} if there exists
a (necessarily unique) thin ideal subGQ containing $x$ and $y$. 
A point $x$ is {\em regular} provided $\{x,y\}$ is regular for every point $y$ which is not collinear with $x$. 

Regularity for lines 
is obtained by interchanging the role of points and lines. 

It is easy to show that if a thick generalized quadrangle of order $(u,v)$ has a regular pair of points (lines), then $v \leq u$ ($u \leq v$).

\subsection{Group notation}

If $\Gamma = (\mP,\mB,\I)$ is a generalized quadrangle, an {\em automorphism} is a permutation of $\mP \cup \mB$ which preserves $\mP$ and $\mB$, and incidence. The set of all automorphisms of $\Gamma$ forms a group under composition of maps, which we denote by $\Aut(\Gamma)$. Each subgroup of $\Aut(\Gamma)$ is called an {\em automorphism group}. 

In general, if $H \leq \texttt{Sym}(\Omega)$ is a group acting on a set $\Omega$, and if $\omega \in \Omega$, then $H_\omega$ denotes the stabilizer of $\omega$ in $H$.

\section{Translation generalized quadrangles}

\subsection{Translation generalized quadrangles}
\label{TGQintro}

A {\em translation generalized quadrangle} (TGQ) is a generalized quadrangle $\Gamma$ for which there  is an abelian automorphism group $A$ that fixes each line incident with some point $x$, while acting sharply transitively on the points not collinear with $x$. We call the latter points the {\em affine points} (with respect to, or relative to, $x$), and $x$ is called {\em translation point}.  The group $A$ is the {\em translation group}.
For the sake of convenience, one should also bear a similar notion in mind for thin quadrangles.

It is possible to show that, as in the case of translation planes, the group $A$ is the unique translation group for the point $x$. There might be more than one translation point, however, but this is not essential for this paper. To make explicit the translation point and/or translation group, we also write $\Gamma^x$ or $(\Gamma,A)$ or $(\Gamma^x,A)$, or even $(\Gamma,x,A)$, for $\Gamma$. So the reader should understand that a morphism 
\begin{equation}
\rho: (\Gamma,x,A) \mapsto (\Gamma',x',A')
\end{equation}
maps $\Gamma$ to $\Gamma'$ and $x$ to $x'$, and necessarily $A$ to $A'$.

We refer the reader to \cite[chapter 8]{PT2} and the monograph \cite{TGQ} for the basics on TGQs.

\subsection{Kantor families}

Keep using the notation of the previous subsection. 

Suppose $z$ is an affine point.
Let 
\begin{equation}
\begin{cases}
\F = \{A_M\ \vert\ M \I z\}; \\
\F^* = \{A_m\ \vert\ x\sim m \sim z\}.
\end{cases}
\end{equation}

Denote $A_m$ also by $A_M^*$. Then for all $L \I z$ we have the following properties:

\begin{itemize}
\item[--(KF1)--]
$A_L \leq A_L^* \ne A_L$;
\item[--(KF2)--]
$A_LA_M^* = A$ for $M \ne L$;
\item[--(KF3)--]
$A_UA_V \cap A_W = \{\he\}$ for distinct lines $U,V,W$;
\item[--(KF4)--]
$\{A^*_L/A_L\} \cup \{A_LA_M/A_L\ \vert\  M \ne L\}$ is a partition of $A/A_L$.
\end{itemize}

Conversely, if families of subgroups $\F$ and $\F^*$ have the properties (KF1)--(KF4) in an abelian group $A$ which a priori need not be related to a TGQ, one says that $(\F,\F^*)$ or $(A,\F,\F^*)$ is a {\em Kantor family} (in $A$). 
Starting from a Kantor family, one can construct  a TGQ, using a natural group coset geometry respresentation \cite{PT2,TGQ}. If one starts from a TGQ as above, it is isomorphic to the reconstructed coset geometry (the isomorphism class of the geometry is independent of the chosen affine point $z$), and the elements of $\F$ and $\F^*$ play the expected role. This representation method was noted for a more general class of generalized quadrangles (namely, for so-called {\em elation generalized quadrangles}) by
Kantor \cite{Ka} in the finite case, and carries over without much change to the infinite case \cite{BaPa}.

\subsection{Symmetry and regularity}

It is easy to see that if $L \I x$, and $M \sim L$, $M$ not incident with $x$, then $A_M$ fixes each line of 
$L^{\perp}$, while acting sharply transitively on the points of $M$ which are not on $L$. One notes that
$A_M$ is independent of the choice of $M \sim L$. Call a line $L$ allowing such a group, in general, an {\em axis of symmetry}. Dually, one defines {\em center of symmetry}. Note that an axis of symmetry is regular, and so is a center of symmetry.

\subsection{Kernel}

The {\em kernel} of the TGQ $(\Gamma,A)$ is the set of endomorphisms $\xi: A \mapsto A$ for which $A_U^{\xi} \leq A_U$ for all $U \I x$, foreseen with natural addition and multiplication (inherited from the endomorphism ring of $A$). We denote the kernel by $\endot(A,A)$.
 
It can be shown that $\endot(A,A)$ is an integral domain \cite{Jos}. For the finite case, one can show more.

\bt[\cite{PT2}, chapter 8]
If $(\Gamma,A)$ is a finite TGQ, its kernel is a commutative field. 
\et

Suppose now that $(\Gamma,A)$ is a TGQ for which $\endot(A,A)$ contains a division ring $k$. Then $A$, and any element of $\F \cup \F^*$, can be seen as vector spaces over $k$. Passing to projective space $\bP$, we obtain a projective representation of $(\Gamma,A)$, and the affine points of $\Gamma$ coincide with the points in an affine space $\bA$ of the same dimension as $\bP$. Denote the hyperplane at infinity by $\bP'$;
the group $A$ acts as the group of translations of $\bP$ with axis $\bP'$. The elements of $\F$ define subspaces of $\bP$, and their intersections with $\bP'$ define a set of subspaces of $\bP'$ which one usually calls ``egg.'' 

For more details we refer to \cite{TGQ}.

\section{Basic structural lemmas}

\medskip
\subsection{Setting}

In this setting, $\mS = \mS^x = (\mS^x,A)$ is a TGQ with translation point $x$ and translation group $A$.  We suppose that the number of points (and so also the number of lines) is not finite. Let $z$ be a fixed affine point, and let $(A,\F,\F^*)$ be the Kantor family in $A$ defined by $z$. 

The kernel $\ell$ is 
defined as above, and we suppose, by way of contradiction, that $k$ is not a division ring.



\subsection{Injectivity}

The proof of the next lemma is essentially the same as in the finite case (which can be found in \cite[chapter 8]{PT2}). We include its proof for the sake of convenience. We will use the representation by coset geometries as in \cite[chapter 8]{PT2}. (It is the only place in the paper where this method is used.)
Below, $\ell^{\times} := \ell \setminus \{0\}$.

\bl
\label{inj}
Each element of $\ell^{\times}$ is injective.
\el

{\em Proof}.\quad
Suppose that $\beta\in \ell$ is such that $a_{0}^\beta = \he$ for some $a_0\in A_0\setminus\{\he\}$, $A_0 \in \F$; then we must show that $\beta = 0$. (The choice of $A_0$ is arbitrary. If $\beta$
has a fixed point not in  $\bigcup_{V \in \F} V$, then it has a fixed point in each $V\setminus\{\he\}$ as well.) Assume the contrary. Choose any element $a_{i}\in A_{i}\setminus\{\he\}$, with $A_i \in \F \setminus \{A_0\}$.
Then the point $a_{0}a_{i}$ is at distance two from $\he$ in the collinearity graph of $\mS^x$. Since $\mS^x$ is thick there exist elements $a_l\in A_l\setminus\{\he\}$ and
$a_k\in A_k\setminus\{\he\}$, where $A_k, A_l \in \F$, $\{A_l,A_k\}\cap \{A_0,A_1\} = \emptyset$ and $A_l\neq A_k$, such that 
\begin{equation}
a_0a_i = a_la_k.
\end{equation}

Letting $\beta$ act yields $a_{i}^\beta=a_{l}^\beta a_{k}^\beta$. First suppose that $a_{l}^\beta=\he$;  then $a_{i}^\beta=a_{k}^\beta$. Since $A_i\cap A_k =\{\he\}$ we obtain that
$a_{i}^\beta = \he$. Analogously $a_k^\beta = \he$ implies that $a_i^\beta = \he$. Next suppose that neither $a_{l}^\beta$ nor $a_{k}^\beta$ equals $\he$. In this case the line $A_la_{k}^\beta$ of $\mS^x$ intersects the line $A_k$ in $a_{k}^\beta\neq \he$
and intersects the line $A_i$ in $a_{i}^\beta\neq \he$. Hence we have found a triangle in $\mS^x$, a contradiction. We conclude that $a_{i}^\beta=\he$, and henceforth that $V^\beta = \he$, for all $V \in \F$. By the
connectedness of $\mS^x$ we know that $A = \Big \langle V\ \vert\ V \in \F \Big\rangle$, and hence it follows that $A^\beta = \he$, that is, 
\begin{equation}
\beta = 0.
\end{equation}  
\eop

\subsection{The GQs $\mS(\alpha,z)$}

Let $\alpha \in \ell^\times$, and $z$ an affine point. To avoid trivialities, suppose that $\alpha$ is not a unit.

It is easy to see that $(A^\alpha,\F^{\alpha},{(\F^*)}^\alpha)$ (obvious notation) defines a Kantor family
in $A^\alpha$. Moreover, the GQ $\mS(\alpha,z)$ defined by this Kantor family is thick, and it is {\em ideal}.
One can derive these properties easily from Lemma \ref{inj}:

\bc
\label{cor4.2}
For each $\alpha \in \ell^{\times}$, $\mS(\alpha,z) \cong \mS$, and $\mS(\alpha,z)$ is an ideal subquadrangle of $\mS$.
\ec
{\em Proof}.\quad
Immediate.
\eop \\

Throughout this paper, we will fix $z$, so that for instance we write $\mS(\alpha)$ instead of $\mS(\alpha,z)$.

\medskip
\section{Adjoining automorphisms}

We keep using the notation of the previous section.

Now for any $n \in \mathbb{N}$, we have that $\mS(\alpha^n) \cong \mS$, and $\mS(\alpha^n) \preceq \mS$. 
Of course, we can also see
$\mS$ as the image under $\alpha$ of a GQ $\mS(\alpha^{-1})$ containing $\mS$. So we have a chain $\mC$ of (embeddings of) isomorphic generalized quadrangles
\begin{equation}
\cdots\ \succeq\  \mS(\alpha^{-1})\ \succeq\ \mS\ \succeq\ \mS(\alpha)\ \succeq\ \mS(\alpha^2)\ \succeq\ \cdots  
\end{equation}

indexed by the integers. Each of these quadrangles contains $x$ and $z$, and each $\mS(\alpha^{u})$ contains $\mS(\alpha^v)$ as an ideal subGQ if $v \geq u$.
Moreover, if $A^{(i)}$ denotes the translation group of $(\mS(\alpha^i),x)$, then $A^{(i)} = A^{\alpha^i}$, and
we have a chain of isomorphic groups 
\begin{equation}
\cdots\ \geq\ A^{(-1)}\ \geq\ A\ \geq\  A^{(1)}\ \geq\ A^{(2)}\ \geq\ \cdots
\end{equation}

It is clear that  
\begin{equation}
\Big(\bigcup_{i \in \mathbb{Z}}\mS(\alpha^i),\bigcup_{i \in \mathbb{Z}}A^{(i)}\Big) =: \Big(\widetilde{\mS(\alpha)},\widetilde{A(\alpha)}\Big)
\end{equation}
 is a translation generalized quadrangle with translation point $x$, and since $\alpha$ stabilizes the chain $\mC$, 
the limit of $\alpha$ is an automorphism of $\widetilde{\mS(\alpha)}$, which we denote by $\widetilde{\alpha}$, that fixes $x$ and $z$ linewise. As such, we have ``adjoined $\alpha$ to $\mS$'' to obtain a new quadrangle for which $\alpha$ has become an automorphism. Note that $\widetilde{\alpha}$ induces a unit in the kernel of $\widetilde{\mS(\alpha)}$ (relative to $z$), and 
so $\widetilde{\alpha} \in \Aut(\widetilde{A(\alpha)})$.

Finally, observe the following.

\bo
For $j \in \mathbb{Z}$, $\widetilde{\alpha}$ induces $\alpha: \mS(\alpha^j) \mapsto \mS(\alpha^{j + 1}) \preceq \mS(\alpha^j)$. Similarly, it induces an element of $\endo^t(A(\alpha^j),A(\alpha^j))$.
\eop
\eo

\br{\rm
As in the introduction, it is clear that $\cup_{i \in \mathbb{Z}}\mS(\alpha^i)$ is isomorphic to the direct limit of the directed system
\begin{equation}
\cdots \xrightarrow{\alpha} \mS \xrightarrow{\alpha} \mS  \xrightarrow{\alpha} \cdots,
\end{equation}
and the similar remark is true for $\cup_{i \in \mathbb{Z}}A^{(i)}$.
}
\er

\medskip
\section{Proof of Theorem \ref{conjrough}}

Now let $\beta$ be a nontrivial element of the kernel of $\widetilde{\mS(\alpha)}$, and suppose it is not a unit. Applying the same construction as above to the pair $(\widetilde{\mS(\alpha)},\beta)$, we obtain a TGQ $\widetilde{\widetilde{\mS(\alpha)}(\beta)}$ by adjoining $\beta$ which we denote by $\widetilde{\mS(\alpha,\beta)}$. It is important to note that the limit of $\alpha$ through this chain still makes sense, and so $\alpha$ ``remains'' an automorphism of $\widetilde{\mS(\alpha,\beta)}$ (see Observation \ref{adjprops}).  
In general, we also use the notation $\widetilde{\mS(\alpha_1,\alpha_2,\ldots)}$, or also $\widetilde{\mS(S)}$ if $S$ is an ordered set.\\

In the next observation, we will use exponential notation for endomorphisms acting on groups and group elements. However if we switch to affine spaces
further on, then usually we use additive notation, so that both $b^m$ and $m\cdot b$ will be used in different contexts ($m \in \mathbb{Z} \setminus \{0\}$) for the same action | this is particularly useful when the endomorphism ``to the power $m$'' can be interpreted as a homothecy (such as in Theorem \ref{TGQembed}).  

Consider the next observation.

\bo
\label{adjprops}
\begin{enumerate}
\item[{\rm (a)}]
If $\alpha \ne \beta$ in $\widetilde{\mS(\alpha,\beta)}$, then they also do not agree on $\widetilde{\mS(\alpha)}$.
\item[{\rm (b)}]
Let $(\Gamma,B)$ be a TGQ with kernel $\wp$, and consider the endomorphism $\widehat{m}: B \mapsto B: b \mapsto b^m$, with $m \in \mathbb{N}$. Then $\widehat{m} \in \wp$. Now let $\ell_m$ be the kernel of $(\widetilde{\mS(\widehat{m})},\widetilde{A(\widehat{m})})$, and let $\epsilon \in \ell_{m}$ be a non-unit. Then 
\begin{equation}
\widehat{m}: \widetilde{A(\widehat{m},\epsilon)} \mapsto \widetilde{A(\widehat{m},\epsilon)}: b \mapsto b^m
\end{equation}
is an element of the kernel $\ell_{m,\epsilon}$ of $\Big(\widetilde{\mS(\widehat{m},\epsilon)},\widetilde{A(\widehat{m},\epsilon)}\Big)$. Moreover, it is a unit, and it induces $\widehat{m}$ on any group $\widetilde{A(\widehat{m})}(\epsilon^j)$, $j \in \Z$. 
In particular, it induces $\widehat{m}$ on $\widetilde{A(\widehat{m})}(\epsilon^0) = \widetilde{A(\widehat{m})}$.
\item[(c)]
The same property as in (b) holds if one replaces $\{ \epsilon \}$ by any finite or countable set $\{ \epsilon_1,\epsilon_2,\ldots \}$.
\end{enumerate}
\eo

{\em Proof}.\quad
Part (a) is obvious.\\

For part (b), first observe that as $B$ is abelian, $\widehat{m}$ indeed is an endomorphism of $B$. Moreover, it maps any subgroup into itself, so it is an element of $\endot(B,B)$. Now consider
\begin{equation}
\widehat{m}: \widetilde{A(\widehat{m},\epsilon)} \mapsto \widetilde{A(\widehat{m},\epsilon)}: b \mapsto b^m.
\end{equation}

We already know it is part of $\ell_{m,\epsilon}$. As each $\widetilde{A(\widehat{m})}(\epsilon^j)$ is a subgroup of $\widetilde{A(\widehat{m},\epsilon)}$, $\widehat{m}$ induces a kernel element of the TGQ $\Big(\widetilde{\mS(\widehat{m})}(\epsilon^j),\widetilde{A(\widehat{m})}(\epsilon^j)\Big)$ for any $j \in \Z$. For each such $j$, we have an isomorphism of TGQs
\begin{equation}
\Big(\widetilde{\mS(\widehat{m})},\widetilde{A(\widehat{m})}\Big) \   \underset{\cong}{\overset{\epsilon^{j}}{\longrightarrow}}  \ \Big(\widetilde{\mS(\widehat{m})}(\epsilon^j),\widetilde{A(\widehat{m})}(\epsilon^j)\Big),
\end{equation}
and it maps $\widetilde{A(\widehat{m})}$ to $\widetilde{A(\widehat{m})}(\epsilon^j) = \epsilon^{-j}\widetilde{A(\widehat{m})}\epsilon^{j}$. For any element $b^{\epsilon^{j}}$ of the latter, we have
\begin{equation}
{(b^{\epsilon^{j}}})^{m} = {(b^m)}^{\epsilon^{j}},
\end{equation}
that is, the following diagram commutes:
\bigskip
\begin{center}
\begin{tikzpicture}[>=angle 90,scale=2.2,text height=1.5ex, text depth=0.25ex]
\node (a0) at (0,3) {$\widetilde{A(\widehat{m})}$};
\node (a1) [right=of a0] {$\widetilde{A(\widehat{m})}(\epsilon^j)$};

\node (b0) [below=of a0] {$\widetilde{A(\widehat{m})}$};
\node (b1) [below=of a1] {$\widetilde{A(\widehat{m})}(\epsilon^j)$};

\draw[->,font=\scriptsize,thick]
(a0) edge node[left] {$\widehat{m}$} (b0)
(a1) edge node[right] {$\widehat{m}$} (b1)
(a0) edge node[auto] {$\epsilon^{j}$} (a1)
(b0) edge node[below] {$\epsilon^{j}$} (b1);

\end{tikzpicture}
\end{center}

as $\epsilon^{j}$ is a group morphism.

It now follows indeed that $\widehat{m}$ indeed is also a unit in $\ell_{m,\epsilon}$.\\

Part (c) follows immediately from part (b).
\eop \\


\bc
\label{stable}
Once a non-unit of type $\widehat{m}$, $m \in \mathbb{N}$, in the kernel $\ell$ of $(\mS,A)$ has been adjoined to $(\mS,A)$, the obtained unit $\widehat{m}$ is stable under adjoints.
\eop \\
\ec



\bo[The quadrangles $\widetilde{\mS(\N)}$]
For any TGQ $(\mS,A)$, the TGQ $(\widetilde{\mS(\N)},\widetilde{A(\N)})$ has a kernel which contains the prime field
$k := \langle 0,1 \rangle$ (which is by definition generated by $0$ and $1$ in the kernel of $(\mS,A)$).
\eo
{\em Proof}.\quad
Each endomorphism $\widehat{m}$ with $0 \ne m \in \N$ is a unit in $\endot(A(\N),A(\N))$.
\eop \\

The field $k$ either is a finite field of prime order, or the field of rationals $\mathbb{Q}$.

\begin{remark}{\rm
Note that $\widetilde{\mS(\mathbb{N})} = \widetilde{\mS(\mathbb{Z})} = \widetilde{\mS(\mathbb{Q})} = \widetilde{\mS(P)}$, with $P$ the set of primes in $\mathbb{N}$.}
\end{remark}

We will only need the first part of the statement of the next observation.

\bo[Transfer of automorphisms]
\label{transfer}
Let $(\Gamma,u,C)$ and $({\Gamma'},{u'},C')$ be TGQs, and suppose $\gamma: (\Gamma,u,C) \mapsto ({\Gamma'},{u'},C')$ is an isomorphism (which maps $u$ to $u'$). Suppose moreover that the kernel of $(\Gamma,u,C)$ contains a prime field $\wp$. Then $\gamma$ induces a semilinear isomorphism between the projective spaces $\bP$ and  $\bP'$ which arise from interpreting $C$, respectively $C'$, as a vector space over $\wp$, respectively $\wp^{\gamma}$. Moreover, this isomorphism maps the egg of $\Gamma$ in $\bP$ to the egg of $\Gamma'$ in $\bP'$.
\eo

{\em Proof}.\quad
The proof is essentially the same as the proof of \cite[Theorem 3.7.2]{TGQ} (which is stated for the case $(\Gamma,u,C) = (\Gamma',u',C')$, and where the finiteness assumption is not used), or as the (concise) proof of \cite[Lemma 1]{BaLuPi}. 
\eop \\

Let $\eta$ be the map | by Observation \ref{transfer}, it is in fact a covariant functor | which associates with a TGQ $(\Gamma,C)$ such as in Observation \ref{transfer} the projective space $\bP$ over $\wp$. Using Observation \ref{transfer} and assuming that some TGQ $(\widetilde{\mS(\mathbb{N})},\widetilde{A(\mathbb{N})})$ is not linear, then with $\gamma$ a non-unit in its kernel, we will show that the commutativity of the diagram
\bigskip
\begin{center}
\begin{tikzpicture}[>=angle 90,scale=2.2,text height=1.5ex, text depth=0.25ex]
\node (a0) at (0,3) {$(\widetilde{\mS(\mathbb{N})},\widetilde{A(\mathbb{N})})$};
\node (a1) [right=of a0] {$\bP$};

\node (b0) [below=of a0] {$(\widetilde{\mS(\mathbb{N})}(\gamma),\widetilde{A(\mathbb{N})}(\gamma))$};
\node (b1) [below=of a1] {$\bP^{\gamma}$};

\draw[->,font=\scriptsize,thick]
(a0) edge node[left] {$\gamma$} (b0)
(a1) edge node[right] {$\gamma$} (b1)
(a0) edge node[auto] {$\eta$} (a1)
(b0) edge node[below] {$\eta$} (b1);

\end{tikzpicture}
\end{center}
(eventually) leads to the desired contradiction in the proof of the next theorem. We first state a lemma which is a variation on Lemma \ref{adjprops}(b).

\bl
\label{adjprops2}
Let $(\mS,A)$ be a TGQ, and suppose its kernel $\ell$ is not a division ring. Consider $(\widetilde{\mS(\N)},\widetilde{A(\N)})$, and let $\widetilde{\ell}$ be the kernel of the latter quadrangle. Let $\beta$ be in $\widetilde{\ell}^\times$, and suppose it is not a unit. Let $r(\beta) \ne 0$ and $h(\beta) \ne 0$ be elements
in $\Z[\beta]$; then both are also elements of $\widetilde{\ell}^\times$ as $\beta \in \widetilde{\ell}$.
\begin{itemize}
\item[{\rm (a)}]
For each $i \in \Z$, $h(\beta)$ can be naturally seen as an element in the kernel of $\widetilde{\mS(\mathbb{N})}(r(\beta)^i)$.
\item[{\rm (b)}]
If $\widetilde{\mS(\mathbb{N})}(r(\beta)^u) \preceq \widetilde{\mS(\mathbb{N})}(r(\beta)^v)$, then $h(\beta)$ for the latter quadrangle
induces $h(\beta)$ for the former. In particular, this is true when $u = 1$. 
\item[{\rm (c)}]
In {\rm (a)} and {\rm (b)}, $\widetilde{\mS(\mathbb{N})}$ can be replaced by $\widetilde{\mS(\mathbb{N},S)}$, with $S$ any subset of $\Z[\beta] \setminus \{0\}$. In particular, {\rm (a)} also holds when  $\widetilde{\mS(\mathbb{N})}(r(\beta)^i)$ is replaced by $\widetilde{\mS(\mathbb{N})(r(\beta))}$.
\end{itemize}
\el
{\em Proof}.\quad
Part (a) follows from the fact that for each $i \in \Z$, the following diagram  (where the horizontal arrows are isomorphisms) commmutes:
\bigskip
\begin{center}
\begin{tikzpicture}[>=angle 90,scale=2.2,text height=1.5ex, text depth=0.25ex]
\node (a0) at (0,3) {$\widetilde{A(\mathbb{N})}$};
\node (a1) [right=of a0] {$\widetilde{A(\mathbb{N})}(r(\beta)^i)$};

\node (b0) [below=of a0] {$\widetilde{A(\mathbb{N})}$};
\node (b1) [below=of a1] {$\widetilde{A(\mathbb{N})}(r(\beta)^i)$};

\draw[->,font=\scriptsize,thick]
(a0) edge node[left] {$h(\beta)$} (b0)
(a1) edge node[right] {$h(\beta)$} (b1)
(a0) edge node[auto] {$r(\beta)^{i}$} (a1)
(b0) edge node[below] {$r(\beta)^{i}$} (b1);

\end{tikzpicture}
\end{center}

as $\Z[\beta]$ is a commutative ring. That is, $h(\beta)$, as an element in the kernel of $\widetilde{\mS(\mathbb{N})}(r(\beta)^i)$, is defined as $r(\beta)^{-i}h(\beta)r(\beta)^i$.

Part (b) follows from the proof of (a).

Part (c) readily follows from the proofs of (a) and (b).
\eop \\

\bt
\label{TGQembed}
Let $(\mS,A)$ be a TGQ, and suppose its kernel $\ell$ is not a division ring. Consider $(\widetilde{\mS(\N)},\widetilde{A(\N)})$.
Then its kernel {\em \underline{is}} a division ring. 
\et
{\em Proof.}\quad
Let $\widetilde{\ell}$ be the kernel of $(\widetilde{\mS(\N)},\widetilde{A(\N)})$.
Suppose $\alpha \in \widetilde{\ell}^{\times}$ is not a unit for $\widetilde{\mS(\N)}$, and consider the sequence
\begin{equation}
\mS \ \hookrightarrow\ \widetilde{\mS(\N)}\ \hookrightarrow\ \widetilde{\mS(\Z[\alpha])}.
\end{equation}

The kernel of $(\widetilde{\mS(\N)},\widetilde{A(\N)})$ contains a prime field $k := \langle 0,1 \rangle$.
We can projectively represent the latter quadrangle in a projective space $\bP$ over $k$. The affine points of the quadrangle correspond to the points of $\bA := \bP \setminus \bP'$, where $\bP'$ is a hyperplane in $\bP$, and $z \in \bA$. 

For any $g(\alpha) \in \Z[\alpha] \setminus \{0\}$ (noting as before that $g(\alpha) \in \widetilde{\ell}^{\times}$, as $\alpha$ is), 
consider the isomorphism
\begin{equation}
g(\alpha): (\widetilde{\mS(\N)},\widetilde{A(\N)}) \longrightarrow (\widetilde{\mS(\N)}(g(\alpha)),\widetilde{A(\N)}(g(\alpha)));
\end{equation}
then we can see the affine points of the latter quadrangle as a subset of $\bA$ which is also an affine space over $k$; also, the stabilizers  in $\widetilde{A(\N)}(g(\alpha))$ of the $\widetilde{\mS(\mathbb{N})}(g(\alpha))$-quadrangle lines incident with the point $z$ define projective $k$-spaces which determine the same egg elements in $\bP'$ as the (same) line stabilizers in $\widetilde{A(\mathbb{N})}$. One observes now (by Observation \ref{transfer}) that $g(\alpha)$ induces a semilinear isomorphism between the projective space $\bP$ and the projective space $\bP^{g(\alpha)} = \widetilde{\bP(\mathbb{N})}(g(\alpha))$ which corresponds to $\widetilde{\mS(\mathbb{N})}(g(\alpha))$. Considering again $\widetilde{\mS(\mathbb{Z}[\alpha])}$, one notes that its kernel contains the commutative field $k(\alpha)$. It has a projective interpretation in a projective space $\widehat{\bP}$ over $k(\alpha)$ 
and we can see $\widehat{\bP}$ naturally as the
union 
\begin{equation}
\bigcup_{f(\alpha) \in \mathbb{Z}[\alpha]}\widetilde{\bP{(\alpha,f(\alpha))}}, 
\end{equation}
where $\widetilde{\bP(\alpha,f(\alpha))}$ is the space corresponding to $\widetilde{\mS(\mathbb{N},\alpha,f(\alpha))}$. (If one orders $\Z[\alpha]$ as $f_0(\alpha), f_1(\alpha),\ldots$, and for each $i \in \mathbb{N}$ we put $F_i := \{f_0(\alpha),\ldots,f_i(\alpha)\}$, then $\widehat{\bP} = \cup_{i \in \mathbb{N}}\widetilde{\bP(F_i)}$, but this amounts to the same thing.)
The ``limit''  $\widehat{\alpha}$ of $\alpha$ is an  automorphism of $\widetilde{\mS(\Z[\alpha])}$ which fixes all TGQ-lines incident with $z$ (and also all the TGQ-lines incident with $x$). One observes that it is induced by a homology of $\widehat{\bP}$ with axis $\widehat{\bP'}$ and center $z$ (this is a general property | see \cite{Jos} | but can readily be deduced from the current proof as well); in the affine space $\widehat{\bA}$ corresponding to the affine points of $\widetilde{\mS(\Z[\alpha])}$, it induces a homothecy with center $z$ and factor ${\alpha}$ (considered as an element of $k(\alpha)$) by Lemma \ref{adjprops2}.

Now fix a $k$-base  $\mB$ 
of $\bA$ such that the coordinates of $z$ are $(0,\ldots,0)$. (Here, a point has coordinates ${(x_i)}_{i \in I}$ for some index set $I$ that contains $1$, and $(0,\ldots,0)$ denotes the point with all entries  $0$.)  By construction, the affine point set of $\widetilde{\mS(\N)}$ consists precisely of all points with coordinates in $k$. The base $\mB$ also is a $k(\alpha)$-base of $\widehat{\bA}$. This follows from the equality $\widehat{\bP} = \cup_{f(\alpha) \in \mathbb{Z}[\alpha]}\widetilde{\bP{(\alpha,f(\alpha))}}$; the affine part of each component $\widetilde{\bP({\alpha,h(\alpha)})}$ in the right-hand side can be expressed  relative to any fixed chosen base of $\bA$ (with coordinate values in the ring $k[\alpha,\frac{1}{\alpha},h(\alpha),\frac{1}{h(\alpha)}]$). In fact, the same equality leads to much stronger consequences, most notably a number of obstructions, since $\widehat{\bP}$ and each $\widetilde{\bP{(\alpha,f(\alpha))}}$ share the same lines on $z$ by construction.
For example, consider a point $z' = {(y_i)}_{I}$ with $y_1 = 1$ and the other entries $0$. Then $\langle z,z'\rangle_{\bP'} \cap \bA$ is a group isomorphic to $(k, +)$, and $\widehat{\alpha}$ induces an endomorphism of the latter (since $\bA^{\widehat{\alpha}} = \bA^{\alpha} \subseteq \bA$ by assumption and $\langle z,z'\rangle_{\bP'}$ is fixed by $\widehat{\alpha}$). Since $k$ is a prime field, it follows that $\widehat{\alpha}$ corresponds to multiplying with a factor from $k^\times$, that is, $\alpha$ already was a unit in the kernel of    $(\widetilde{\mS(\N)},\widetilde{A(\N)})$ (and even contained in $k$), contradiction.

The theorem is proved. 
\eop \\

\bc
\label{corembed}
If the kernel of a TGQ $(\mS,A)$ {\em contains} a division ring, it {\em \underline{is}} a division ring.
\ec

{\em Proof.}\quad
Let $(\mS,A)$ be as in the statement, and let $\ell$ be its kernel. If $\wp$ is a division ring contained in $\ell$, it follows that 
$(\mS,A) = (\widetilde{\mS(\N)},\widetilde{A(\N)})$ (since $\langle 0,1 \rangle \leq \wp$ is a field), so it suffices to apply Theorem \ref{TGQembed}.  
\eop \\

This ends the proof of Theorem \ref{conjrough}.

\br
{\rm
Theorem \ref{TGQembed} can also be obtained in an entirely different fashion, still assuming the existence of $\alpha$ by way of contradiction, 
by using \S\S \ref{subs9.1} (the paragraph after Observation \ref{obs9.1}).
}
\er

\medskip
\section{The characteristic of a TGQ}

Let $(\mS,A)$ be a TGQ, and let $(\widetilde{\mS},\widetilde{A})$ be an ambient TGQ. Then we call 
the characteristic of the prime field of its kernel  $k$ the {\em characteristic} of $(\mS,A)$. We need to show that this notion is independent of the choice of ambient TGQ.

Let the characteristic be $0$; then the endomorphism $1$ group theoretically generates $(\Z,+)$ in $\endot(A,A)$. If $m \in \Z$ is not a unit in $\endot(A,A)$, then we have a chain of strict embeddings
\begin{equation}
\cdots\ \succeq\  \mS(1/\widehat{m})\ \succeq\ \mS\ \succeq\ \mS(\widehat{m})\ \succeq\ \mS(\widehat{m^2})\ \succeq\ \cdots  
\end{equation}

We will see in the next section that this implies that any other ambient TGQ also yields characteristic $0$. And we will observe that 
if the characteristic is $p > 0$ for some ambient TGQ, it is so for {\em every} ambient TGQ. (Alternatively, it is easy to see that ambient TGQs are in fact unique by construction.)

\medskip
\section{TGQs in positive characteristic}

Now let $(\mS,A)$ be a TGQ with kernel $\ell$, and suppose that $0$ and $1$ generate $\langle 0,1 \rangle \cong \bF_p$, $p$ a prime, in the kernel $k$ of $(\widetilde{\mS(\N)},\widetilde{A(\N)})$. Represent the latter quadrangle in a projective space $\bP$ over $\bF_p$, with affine point set $\bA = \bP \setminus \bP'$, with $\bP'$ a hyperplane in $\bP$. Suppose $\lambda \in \Z \setminus \{0\}$, and interpret $\lambda \in \bF_p^{\times}$ through the ring morphism $\Z \mapsto \bF_p$; 
suppose it is not a unit in $\ell$ (and note that it is contained in $\endot(A,A)$). Then $\mS^{\lambda^k} \ne \mS^{\lambda^j} \preceq \mS^{\lambda^k}$ for all positive $j, k \ne 0$ with $k < j$. This contradicts the fact that $\lambda^{p - 1} = \id$. 

It follows that $\mS = \widetilde{\mS(\N)}$, and by Corollary \ref{corembed} we conclude that the kernel of $(\mS,A)$ is a division ring.

This solves Conjecture \ref{conjlinear} in positive characteristic.

\medskip
\section{Remarks on characteristic $0$}

We now suppose that $(\mS^x,A)$ is a TGQ in characteristic $0$. As before, $z$ is a fixed affine point.

\subsection{The core of $(\mS^x,A)$}
\label{subs9.1}

Define the {\em core} $\Delta$ of $(\mS^x,A)$ to be  the intersection of all ideal subGQs (both thick and thin) of $\mS$, which contain $x$ and $z$. 

\bo
\label{obs9.1}
Suppose $\Delta$ is thick.
Then $\Delta$ is a TGQ with translation point $x$ which is linear. 
\eo

{\em Proof}.\quad
Let $\kappa \ne 0$ be an element of the kernel of $\Delta$. Then if $\Delta^\kappa \ne \Delta$, we have that $\Delta^\kappa$ is a proper ideal subGQ
of $\Delta$ that contains $x$ and $z$, contradiction.
\eop \\

Now $\Delta$ also is an ideal subquadrangle of $(\widetilde{\mS},\widetilde{A})$ containing $x$ and $z$. Represent the latter in a projective space $\bP$ over the prime field $\mathbb{Q} \cong \langle 0,1\rangle \leq k$.
If $\Delta$ is thick, then all projective lines on $z$ which contain at least one other affine point of $\Delta$  must be completely contained in $\Delta$, by the action of $\langle 0, 1\rangle$ by homothecies. In fact, the core precisely singles out all projective lines over $\langle 0,1 \rangle$ which are incident with $z$ and completely contained in $\mS$ (and such lines always exist). 

Let $A_{\Delta} \leq A$ be the translation group of $\Delta$ relative to the translation point $x$. Then $A$ forms a vector space over $\langle 0,1\rangle$. Denote the corresponding affine space by $\bA_{\Delta}$.

\begin{conjecture}\label{conjend}{\rm
$\dim{\bA_{\Delta}} = \dim{\bA}$.
}\end{conjecture}

If Conjecture \ref{conjend} is true, it immediately follows that 

\begin{equation}
\Delta = \mS = \widetilde{\mS}.
\end{equation}

So if Conjecture \ref{conjend} is true, $\Delta$ must be a thin ideal subGQ, implying that $x$ is a regular point in $\mS$. Since $\mS^x$ is a TGQ, $x$ is a center of symmetry and 
its symmetries are contained in $A$. 


\subsection{The finite case}

In the finite case, this combinatorial situation | a center of symmetry  which is incident only with axes of symmetry | 
already would lead to the desired contradiction: if a finite thick generalized quadrangle of order $(u,v)$ has a regular point which is incident with only regular lines,
then $u = v$ and $u$ is even | see \cite[1.5.2]{PT2}. For (infinite) Moufang quadrangles, the existence of a regular point and a regular line also is sufficient to conclude that one is working in characteristic $2$, see Tits \cite{QDM} and also \cite{TiWe}. But with only local Moufang conditions, this is no longer true. For, consider any infinite field $k$; then in $\bP^2(k)$ there 
exists a set $S$ of points such that any line meets $S$ in precisely two distinct elements. Embed $\bP^2(k)$ as a hyperplane in $\bP^3(k)$, choose a point 
$z \in \bP^3(k) \setminus \bP^2(k)$, and choose a point $c \in S$. For any $s \in S$, define the group $T_s := T_{zs}$, where $T$ is the translation group of the affine space $\bP^3(k) \setminus \bP^2(k)$. Then $\F_S = \{ T_s\ \vert\ s \in S \setminus \{s\}\}$ defines a Kantor family in $T$, where $\F^* = \{ T^*_s\ \vert\ T^*_s = T_s \oplus T_c, s \ne c\}$. The reader notices that from this Kantor family a TGQ arises with kernel $k$, and regular translation point.

\subsection{Restricting non-linear TGQs in characteristic $0$}

Finally, we observe the following result.

\bo
If $(\mS^x,A)$ is a TGQ which is not linear, then we have the following properties.
\begin{itemize}
\item[{\rm (CH)}]
It lives in characteristic $0$.
\item[{\rm (CO)}]
It properly and ideally contains a linear TGQ $(\Delta,\underline{A})$ (the core).
\item[{\rm (AM)}]
It is properly and ideally contained in a linear TGQ $(\widetilde{\mS},\widetilde{A})$ (the ambient TGQ). 
\item[{\rm (DI)}]
$[\widetilde{A} : A] = \vert \mathbb{N} \vert$.
\end{itemize}
If $(u,v)$ is the order of $\mS$, then $u = v \not\in \mathbb{N}$; this is also the order of the ambient TGQ $\widetilde{\mS}$.
\eo

{\em Proof}.\quad
We only need  to prove the last statement. Since $\mS^x$ is a TGQ, it contains regular lines, and so $v \geq u$ (see \S \ref{reg}). If $\Delta$ is not thick, we know $x$ is a regular point of $\mS$, so that $v \leq u$. Next, suppose that $\Delta$ is thick. 
Now let $X$ be a  line of $\widetilde{\mS}$ which is not contained in $\mS$. (Note that $\widetilde{\mS} \ne \mS$ since $\mS$ is not linear.)
Project every point of $\mS$ on $X$; we obtain a set $S$ of lines of $\mS$ which partitions the point set of $\mS$, and which all meet $X$. 
Since the number of points of $\mS$ is $(u + 1)(uv + 1) = v$, the order of $S$ is $v$. Now observe that by construction of $\widetilde{\mS}$,  lines of $\mS$ and $\widetilde{\mS}$ have the same number $u$ of points  (since each line of $\widetilde{\mS}$ is a countable union of lines of GQs isomorphic to $\mS$). It follows that 
\begin{equation}
u = \vert X \vert\ \geq\ \vert S \vert = v. 
\end{equation} 
We conclude that $u = v$. 
\eop \\


\subsection{Final remarks}

If $\Delta$ is thick, it is not hard to prove that $\Delta$ and $\widetilde{\mS}$ have the same kernel. 

Also, note that by Corollary \ref{corembed}, we have that a TGQ $\mS$ with kernel $\ell$ is nonlinear if and only if some $\widehat{m}$ in $\ell^\times$ ($m \in \Z$) is not a unit. Since we know that $\ell$ has to have characteristic $0$, we conclude that for solving Conjecture \ref{conjlinear},
it is sufficient to handle the cases where $\ell$ is a subring of $\mathbb{Q}$.

\end{document}